\documentclass[reqno]{amsart}

\usepackage{amsmath, amsfonts, amsthm, amssymb, color,  graphicx, mathrsfs, cite}
\usepackage{comment,stmaryrd}

\theoremstyle{plain}
\newtheorem{theorem}{Theorem}[section]
\newtheorem{lemma}{Lemma}[section]
\newtheorem{proposition}{Proposition}[section]
\newtheorem{corollary}{Corollary}[section]

\theoremstyle{definition}
\newtheorem{definition}{Definition}[section]

\theoremstyle{remark}

\newtheorem{example}{Example}[section]

\numberwithin{equation}{section}
\allowdisplaybreaks

\usepackage{ifpdf}
\ifpdf \usepackage[colorlinks=true, citecolor=blue, linkcolor=blue, urlcolor=blue]{hyperref} \fi

\def\thrm{\begin{theorem}}
\def\thrml#1{\begin{theorem}\label{#1}}
\def\ethrm{\end{theorem}}
\def\lmm{\begin{lemma}}
\def\lmml#1{\begin{lemma}\label{#1}}
\def\elmm{\end{lemma}}
\def\dfntn{\begin{definition}}
\def\dfntnl#1{\begin{definition}\label{#1}}
\def\edfntn{\end{definition}}
\def\crllr{\begin{corollary}}
\def\crllrl#1{\begin{corollary}\label{#1}}
\def\ecrllr{\end{corollary}}
\def\xmpl{\begin{example}}
\def\xmpll#1{\begin{example}\label{#1}}
\def\exmpl{\end{example}}
\def\nmrt{\begin{enumerate}}
\def\enmrt{\end{enumerate}}
\def\qtn{\begin{equation}}
\def\qtnl#1{\begin{equation}\label{#1}}
\def\eqtn{\end{equation}}
\def\prpstn{\begin{proposition}}
\def\prpstnl#1{\begin{proposition}\label{#1}}
\def\eprpstn{\end{proposition}}

\def\tm#1{\item[{\rm (#1)}]}
\def\proof{{\bf Proof}.\ }
\def\eprf{\hfill$\square$}

\DeclareMathOperator{\supp}{supp}

\newcommand{\F}{\mathbb{F}}

\newcommand{\Z}{\mathbb{Z}}

\def\qaq{\quad\text{and}\quad}

\def\lg{\langle}
\def\rg{\rangle}

\def\cZ{\mathcal {Z}}

\def\cA{\mathcal {A}}

\def\cD{{\mathcal D}}

\def\qaq{\quad\text{and}\quad}

\begin{document}

\title{Schur rings over infinite dihedral group }

\author{Gang chen}
\address{School of Mathematics and Statistics, Central China Normal University, Wuhan 430079, China}
\email{chengangmath@mail.ccnu.edu.cn}

\author{Jiawei He}
\address{School of Mathematics and Information Sciences, Nanchang HangKong University, Nanchang 330063, China}
\email{hjwywh@mails.ccnu.edu.cn}

\author{Zhiman Wu}
\address{Fengcheng No.9 Middle School, Yichun 331100, China}
\email{1793938418@qq.com}

\begin{abstract}

Schur rings over the infinite dihedral group $\cZ\rtimes\mathcal{Z}_2$ are studied according to properties of Schur rings over infinite groups and the classification of Schur rings over infinite cyclic groups.  Schur rings over $\mathcal{Z}\rtimes{\mathcal{Z}}_2$ are  classified under the assumption that $\cZ$ is an $\cA$-subgroup. Those Schur rings are proved to be traditional.
\end{abstract}


\date{}

\maketitle


\section{Introduction}\label{in}

Schur rings are a type of subrings of group rings which are  determined by a partition of the underline group, and they were originally developed by Schur \cite{Sch} and Wielandt \cite{Wie} in order to study permutation groups as an alternative to character theory. In addition, theory of Schur rings  has  been used in many other applications. For the applications of Schur rings to association
schemes, see for example, subsections 2.4 and 4.4 of \cite{CP}.

\medskip

Classification of Schur rings over a given group  $G$ has been a subject of much research interest.
Most papers have studied Schur rings over finite groups.
Leung and Man have classified Schur rings over finite cyclic groups. There are four types of Schur rings over finite cyclic groups which are: partitions induced by direct product decomposition of groups (Schur rings of direct produt type);  partitions induced by orbits of automorphism subgroups (orbit Schur rings); partitions induced by cosets of normal subgroups (Schur rings of wedge product type); trivial Schur ring. We say that a Schur ring over a group is a {\it traditional Schur ring} if it is one of these four types.

\medskip

Some scholars have extended some properties of Schur rings over finite groups to Schur rings over infinite groups and have classified Schur rings over infinite cyclic groups,  the group $\mathcal{Z}\times {\mathcal{Z}_2}$  $(\mathcal{Z}$ is an infinite cyclic group, $\mathcal{Z}_2$ is a cyclic group of order 2)~and torsion-free locally  cyclic groups. In  particular, the classification of Schur rings over infinite cyclic groups and the main
results of Leung and Man tell us that all Schur rings over cyclic groups are traditional.

\medskip
 Although general  structure theorems are available for Schur
rings over finite cyclic groups, the problem of classification of Schur rings over infinite groups is proved to be a difficult problem.  Some special cases have been considered, for example, Bastian et al.  \cite{BBHMT} study the Schur rings over some types of infinite groups. According to their results,  there are nontraditional Schur
rings over free groups and free products. Recently, Chen et al. \cite{CH} classify
 Schur rings over $\mathcal{Z} \times \mathcal{Z}_3$. In particular, all these Schur rings are proved to be traditional. In this paper we consider the  classification of Schur rings over $\mathcal{Z}\rtimes {\mathcal{Z}_2}$.

\medskip
Throughout let $G$ denote an arbitrary group with identity element $1$, and $\mathbb{F}$ an arbitrary field of characteristic $0$ and we assume that it contains $\Z$ as a subring. The group algebra over $G$ with coefficients in $\F$ is denoted $\F[G]$. The infinite dihedral group $\langle z,s~:~s^2=1,s^{-1}zs=z^{-1}\rangle$ is denoted by $\cD_{\infty}$.

\medskip

As we shall see, the main results of the paper tells us that when $\langle z\rangle$ is an $\cA$-subgroup for a Schur ring over $\cD_{\infty}$. Then can get the exact classification of $\cA$, which are prove to be traditional. The following results implies that under some condition,  $\langle z\rangle$ is really an $\cA$-subgroup.

\medskip

\prpstnl{1136a} Let $\cA$ be a Schur ring over $\cD_{\infty}$, $\cD$  the basic set of $\cA$, and $C$ the basic set containsing $z$. Then the following statements hold:
\nmrt
\tm{i} If  $C$ is contained in $\cZ$, then $\langle z\rangle$ is an $\cA$-subgroup.

\tm{ii} Suppose that $C$ contains elements of the form $z^is$. Then there exist at least two elements $z^is$ and $z^js$ with $i\ne j$ in $C$.

\tm{iii} Set

\enmrt
$$
l:=max_{D\in\cD}\{\:i-j\::\:z^is,z^js\in D\:\}
~~
and
~~
t:=min_{D\in\cD}\{\:|i-j|\::\:z^is,z^js\in D\:\}.
$$
Suppose that  ~$2\nmid l$ and~$t>2$, then $\cZ$ is an $\cA$-subgroup.
\eprpstn

The following are the main results of the present paper.

\thrml{1451b} Let $\cA$ be a Schur ring over~$\cD_{\infty}$. Assume that $\cZ=\langle z\rangle$~is an~$\cA$-subgroup.  Then $\cA$ is isomorphic to one of the following:
\nmrt
\tm{i} $\cA=\F[\cZ\rtimes{\cZ}_2];$

\tm{ii} $\cA=\mathrm{Span}_\F\left\{\:{\underline C\::\:C \in \cD\:} \right\}$,
where $\cD=\{\{z^js,z^{i-j}s\}:\:j
\in \Z\:\}\cup \{\{z^j,z^{-j}\}:j\in \Z\}$~for some~$i\in \Z\backslash \{0\}$.
\enmrt
In addition the Schur ring in $(ii)$ is an orbit Schur ring.
\ethrm

\thrml{2305b} Let
$$
\cA=\mathrm{Span}_{\mathbb{F}}\{~\underline{C}~:~C\in \mathcal{D}~\},
$$
where
\begin{displaymath}
	\cD=\{\{z^{i},z^{-i},z^is,z^{-i}s\}:~2\nmid i,i\in \mathbb{Z}~\}\cup\{\{z^is^j, z^{-i}s^j\}:2\mid i,~i\in \mathbb{Z},j\in\{0,1\}\}.
\end{displaymath}
Then  $\cA$ is a Schur ring over $\cD_{\infty}$ such that $\langle z\rangle$ is not an $\cA$-subgroup and furthermore $\cA$ is not traditional.

\ethrm

\section{Preliminary}\label{Pr}
For any finite non-empty subset $C\subseteq G$, denote $\underline{C}=\sum_{g\in C}g$ and call $\underline{C}$ a {\it simple quantity}. In addition, we define $C^*=\{g^{-1}: g\in C\}$. A partition $\cD$ of $G$ is said to be  {\it finite support} provided any $C\in \cD$ is a finite subset of $G$. Our notation for Schur ring is taken from \cite{Po}.

 \subsection{Schur Rings}
\dfntnl{1055a} Let $\cD$ be a finite support partition of $G$ and $\cA$ a subspace of $\mathbb{F}[G]$ spanned by $\{\underline{C}: C\in \cD\}$. We say that $\cA$ is a {\it Schur ring} (or $S$-ring) over $G$ if
 \nmrt
  \tm{1} $\{1\} \in \cD$,
  \tm{2} for any $C\in \cD$, $C^*\in \cD$,
  \tm{3} for all $C, D\in \cD$, $\underline{C}\cdot\underline{D}=\sum_{E\in \cD}\lambda_{CDE}\underline{E}$, where all but finitely many $\lambda_{CDE}$ equal $0$.
 \enmrt
\edfntn

For a Schur ring $\cA$ over $G$, the associated partition $\cD$ is denoted $\cD(\cA)$ and each element in $\cD(\cA)$ is called a {\it basic  set}
of $\cA$.  Each $\lambda_{CDE}$ is called a {\it structure constant}. We say that a subset $C$ of $G$ is an {\it $\cA$-set} if $C$ is a union of some basic sets of $\cA$. An {\it $\cA$-subgroup} is simultaneously an $\cA$-set and a subgroup.

\medskip

Suppose $C=\sum_{g\in G}\alpha_gg$ and $D=\sum_{g\in G}\beta_gg\in \F[G]$ (note that here only finitely many nonzero coefficients $\alpha_g$ and $\beta_g$). Then define
$C^*=\sum_{g\in G}\alpha_gg^{-1}$ and $\supp(C)=\{g: \alpha_g\ne 0 \}$. Additionally,  the {\it Hadamard product} is defined as $C\circ D=\sum_{g\in G}\alpha_g\beta_gg$. \medskip

\lmml{1718c} Let $\cA$ be a Schur ring over a group $G$. Suppose $C=\{g\}$ is a basic set. Then for any basic set $D$, $gD$ is a basic set.

\elmm

The following theorem was proved by Wielandt in \cite{Wie2} when $G$ is finite. The general case appeared as \cite[Corollary 2.12]{BBHMT}.

\thrml{1505b} Let $G$ be a group and $\cA$ a subring of $\mathbb{F}[G]$ which is spanned by $\{\underline{C}: C\in \cD(\cA)\}$ as an $\F$-vector space. Then $\cA$ is a Schur ring over $G$ with basic set $\cD(\cA)$ if and only if $\cA$ is closed under $\circ$ and $*$, $1\in \mathcal{D}(\cA)$, and for all $g\in G$ there exists some $C\in \mathcal{D}(\cA)$ such that $g\in \supp(C)$.
\ethrm

Let $f: \F\rightarrow \F$ be a function such that $f(0)=0$. For any $\alpha=\sum_{g\in G}\alpha_gg$, we set
$$
f(\alpha)=\sum_{g\in G}f(\alpha_g)g.
$$

The following proposition was proved by Wielandt in \cite[Proposition 22.3]{Wie2}. The general case appeared as \cite[Propostion 2.4]{BBHMT}.

\prpstnl{913a}Let $\cA$ be a Schur ring over a group $G$ and the function be as above. Then
$f(\alpha)\in \cA$ whenever $\alpha\in \cA$.
\eprpstn

Let $f$ take value $1$ at one non-zero number and $0$ at other numbers.  We get the following statement, which is known as the Schur-Wielandt principle; see \cite[Corollary 1.10]{Po}.

\crllrl{10131c} Let $\cA$ be a Schur ring over $G$. For any $\alpha\in \cA$, suppose $g\in \supp({\alpha})$ with $\alpha_g=c\ne 0$. The set
$$
\{g\in G: \alpha_g=c\}
$$
is an $\cA$-set.
\ecrllr

The following two propositions, which were first proved by Wielandt in \cite{Wie2} and were generalized in \cite{BBHMT},  tell us how to generate $\cA$-subgroups in any Schur ring $\cA$.

\prpstnl{956b} Let $\cA$ be a Schur ring over $G$. Let $\alpha\in \cA$ and
$$
{\rm Stab}(\alpha)=\{g\in G: g \alpha=\alpha\}.
$$
Then ${\rm Stab}(\alpha)$ is an $\cA$-subgroup.
\eprpstn

\crllrl{1539b} Let $\cA$ be a Schur ring over an infinite group $G$. Then for any basic set $C$
$$
{\rm Stab}(C)=\{g\in G: g C=C\}
$$
is a finite $\cA$-subgroup.
\ecrllr
\proof Observe that $gC=C$ if and only if $g\underline{C}=\underline{C}$. It then follows from Proposition \ref{956b} that ${\rm Stab}(C)$ is an $\cA$-subgroup. Since we are assuming that $C$ is a finite set, ${\rm Stab}(C)$ is a finite subgroup of $G$.
\eprf

\prpstnl{956b} Let $\cA$ be a  Schur ring over $G$. Let $\alpha\in \cA$ and $H=\lg \supp(\alpha)\rg$. Then~$H$ is an $\cA$-subgroup.

\eprpstn

Suppose that $H$ is an $\cA$-subgroup, then $\{C: C\in \cD(\cA), C\subseteq H\}$
consists of  a basic set of a Schur ring over $H$. It is denoted $\cA_H$.

\medskip

When we have normal $\cA$-subgroup, we can construct  Schur ring over the factor group as shown in the following lemma.

 \lmml{1440a}(\cite[Lemma 1.2]{LM}) Let $\varphi: G\rightarrow H$ be a group homomorphism with ${\rm Ker}(\varphi)=K$ and $\cA$ be a Schur ring over $G$. Suppose that $K$ is an $\cA$-subgroup. Then the image $\varphi(\cA)$ is a Schur ring over $\varphi(G)$ where $\cD(\varphi(\cA))=\{\varphi(C): C\in \cD(\cA)\}$.
\elmm

In particular, if $H=G/K$, the factor group of $G$ over $K$, and $\varphi$ is the natural homomorphism,  the corresponding $S$-ring is denoted $\cA_{G/K}$.

\medskip

\thrml{1514a}(\cite[Theorem  3.3]{BBHMT}) The only Schur rings over the infinite cyclic group $\langle z \rangle$ are either discrete or symmetric Schur ring. Here each basic set of the discrete Schur ring is a singleton set and each basic set of the  the symmetric Schur rings has the form $\{z^i, z^{-i}\}$ for some nonnegative integer $i$.
\ethrm

\subsection{Traditional Schur rings}
\
\medskip

\textbf{Discrete and trivial Schur rings}.
The group ring $\mathbb{F}G$ itself forms a Schur ring over $G$, which is called
the {\it discrete Schur ring} over $G$. When $G$ is a finite group, the trivial partition $\{1,G \}$ produces
a Schur ring which is known as the {\it trivial Schur ring} over $G$.

\medskip
\textbf{Tensor products}. Assume that $G =H \times K $, $\mathcal{A}_H$ and $\mathcal{A}_K$ are Schur rings over $H$ and $K$, respectively. Moreover, the
set $$
\{CD : C\in \mathcal{D}(\cA_H), D\in \mathcal{D}(\cA_H)\}
$$
consists of the basic sets of a Schur ring over $G$, which is first proved by Wielandt in \cite{Wie}. The corresponding Schur ring is called the {\it tensor product} of $\mathcal{A}_H$ and $\mathcal{A}_K$.

\medskip

\textbf{Orbit Schur rings}. If $\mathcal{K}$ is a finite subgroup of $\operatorname{Aut}(G)$, the set of elements of $\mathbb{F}[G]$ fixed by $\mathcal{K}$ is a Schur ring over G, denoted $\mathbb{F}[G]^{\mathcal{K}}$ and called the {\it orbit Schur ring} associated with $\mathcal{K}$. Obviously, the group ring $\mathbb{F}G$ itself also forms an orbit Schur ring.

\medskip

\textbf{Wedge products}. A Schur ring $\mathcal{A}$ over $G$ is a {\it wedge product} if there exist nontrivial proper $\mathcal{A}$-subgroups $H, K$ such that $K \leq H, K \unlhd G$, and every basic set outside $H$ is a union of $K$-cosets. In this case, the series
$$
1<K \leq H<G
$$
is called a wedge-decomposition of $\mathcal{A}$. Furthermore, we write $\mathcal{A}=\mathcal{A}_1 \wedge \mathcal{A}_2$, where $\mathcal{A}_1=\mathcal{A}_H$ and $\mathcal{A}_2=\mathcal{A}_{G / K}$.

\medskip

A Schur ring over a group $G$ is called {\it traditional} if it is either a trivial Schur ring (when $G$ is finite), or an orbit Schur ring, or a tensor product of Schur rings over smaller subgroups, or a wedge product. A group $G$ is said to be {\it traditional} if each Schur ring over it is traditional.

\medskip

\section{Proof of the main Results}

{\bf Proof of Proposition \ref{1136a}}\,
If~$z^is\notin C$~for all $i\in \mathbb{Z}$, then~$\cZ=\langle C\rangle$~is an~$\cA$-subgroup. Then statement $(i)$ follows.

\medskip

 Next, assume that ~$z^is\in C$~for some $i\in \mathbb{Z}$. If there is only a unique $0\neq i\in \mathbb{Z}$ such that $z^is$ is contained in $C$. Observe that~$z^is\in C\cap C^{*}$, and we  conclude that $C=C^{*}$. This yields that
$$
\underline{C\backslash \{z^is\}}\cdot z^is=z^is\cdot\underline{C\backslash \{z^is\}}.
$$
Also, for any odd prime number $p$,
\begin{displaymath}
	\begin{split}
		\underline{C^{p}}&={(\underline{C\backslash \{z^is\}}+z^is)}{(\underline{C\backslash \{z^is\}}+z^is)}...{(\underline{C\backslash \{z^is\}}+z^is)}\\
		&\equiv ({\underline{C\backslash \{z^is\}}})^{p}+(z^is)^{p}\\
		&\equiv ({\underline{C\backslash \{z^is\}}})^{(p)}+(z^is)^{p}\\
		&={\underline{C}^{(p)}}~~~~~~~~~~~~~~~~~~~~~~~~~(mod~~p),
	\end{split}
\end{displaymath}
and so~$C^{(p)}$~is an~$\cA$-set. This means that $C\subseteq C^{(p)}$ as
${(z^is)}^{p}=z^is$, and hence $C=C^{(p)}$ because
$$
|C|\leq |C^{(p)}|\leq |C|<\infty.
$$
This implies that $z^{p}\in C$ for any odd prime number $p$, which contradicts the finiteness of $C$.
Thus there exists a positive integer $k$ greater than or equal to 2 such that
$$
z^{m_1}s,z^{m_2}s,...,z^{m_k}s\in C,
$$
~where~$m_1,m_2,...,m_k\in \mathbb{Z}$~and~$m_1<m_2<...<m_k.$ Thus, statements $(ii)$ follows.

\medskip

Now, we will prove conclusion $(iii)$.
Note that
\begin{displaymath}
	\begin{split}
		&\underline{C\backslash \{z^{m_1}s+z^{m_2}s+...+  z^{m_k}s\}}\cdot(z^{m_1}s+z^{m_2}s+...+z^{m_k}s)\\=
		&(z^{m_1}s+z^{m_2}s+...+z^{m_k}s)\cdot\underline{C\backslash \{z^{m_1}s+z^{m_2}s+...+z^{m_k}s\}}.
	\end{split}
\end{displaymath}
Thus,
\begin{displaymath}
	\begin{split}
		{\underline{C}}^2&=\{\underline{C\backslash \{z^{m_1}s+z^{m_2}s+...+  z^{m_k}s\}}+(z^{m_1}s+z^{m_2}s+...+  z^{m_k}s)\}^2\\
		&\equiv (\underline{C\backslash \{z^{m_1}s+z^{m_2}s+...+  z^{m_k}s\}})^2+(z^{m_1}s+z^{m_2}s+...+  z^{m_k}s)^2\\
		&\equiv ({\underline{C\backslash \{z^{m_1}s+z^{m_2}s+...+z^{m_k}s\}}})^{(2)}+(z^{m_1}s+z^{m_2}s+...+  z^{m_k}s)^2~~~~~~~(mod~~~2).
	\end{split}
\end{displaymath}
Let $K$ denote the set of all elements of~${\underline{C}}^{2}$~with odd coefficients. Clearly, $K\subseteq \mathcal{Z}$ and is an $\cA$-set.
Observe that
\begin{displaymath}
	\begin{split}
		&(z^{m_1}s+z^{m_2}s+...+z^{m_k}s)^2\\
		=&(1+z^{m_1-m_2}+z^{m_1-m_3}+...+z^{m_1-m_k})+\\
		&(1+z^{m_2-m_1}+z^{m_2-m_3}+...+z^{m_2-m_k})+...+\\
		&(1+z^{m_k-m_1}+z^{m_k-m_2}+...+z^{m_k-m_{k-1}}),
	\end{split}
\end{displaymath}
by the assumption in statement $(iii)$, one can see that $m_k-m_1$ is odd and
$$
m_{i+1}-m_i> 2, i=1, \ldots, k-1.
$$
These facts yield that
$$
\{z^{m_1-m_k},z^{m_k-m_1}\}\subseteq K,~ \{z^2,z^{-2}\}\subseteq K.
$$
Therefore~$\cZ=\langle K\rangle$~is an $\cA$-set, and hence we complete the proof of Proposition~\ref{1136a}.~\eprf

\medskip

{\bf Proof of Theorem \ref{1451b}}\,  In the sequel, we fix the following notation:
$$
G=\cZ \rtimes \cZ_{2}=\lg z\rg \rtimes\lg s\rg,
$$
and $\cA$ is  a Schur ring over $G$ such that $\cZ$ is an $\cA$-subgroup.

\medskip

By Theorem \ref{1514a},~we consider the following two:

 \medskip

{\bf Case 1.}\, Suppose ~$\cA_{\cZ}=\F[\cZ]$ is a discrete Schur ring.

\medskip

 We claim that each basic set in $\cD(\cA)$ is a singleton set.  Otherwise there exists~$i\neq j\in \Z$ such that $z^is,z^js\in C$ for some basic set $C$. It then follows that $z^{i-j}C$~ is a basic set by Lemma \ref{1718c}. Observe that
$$
z^{i}s=z^{i-j}\cdot z^{j}s\in z^{i-j}C,
$$
and so~$z^{i-j}C=C$. This yields that $z^{i-j}$ belongs to ${\rm Stab}(C)$ and which would contain the infinite subgroup $\langle z^{i-j} \rangle$. This is a contradiction to Lemma \ref{1539b}.

\medskip

As a consequence, $\cA=\F[G]$.

\medskip

{\bf Case 2.}\, Suppose~$\cA_{\cZ}={\F[\cZ]}^{\pm}$ is the symmetric Schur ring. Let  $D$ be the basic set containing~$s$.

 \medskip

{\bf Subcase 2.1.}\, Suppose $|D|=1$. It then follows that
$$
\cD(\cA)=\{\{z^js,z^{-j}s\}:j\in \Z\}\cup \{\{z^j,z^{-j}\}|j\in \Z\}.
$$
Thus case (ii) in Theorem \ref{1451b} occurs with $i=0$.

\medskip

{\bf Subcase 2.2.}\, Suppose $|D|=2$. Then there exists~$i\in \Z\backslash \{0\}$~such that $D=\{s,z^is\}$.

\medskip

{\bf Claim 1.}\, There exists a basic set of the form $\{z^{i_0}s\}$ for some $i_0\in \mathbb{Z}\backslash \{0\}$ if and only if $i$ is even. Furthermore, if this happens then
$\{z^{\frac{i}{2}}s\}$ is the unique basic set of size one and
$$
\cD=\{\{z^{j+\frac{i}{2}}s,z^{\frac{i}{2}-j}s\}:j\in \Z\}\cup \{\{z^j,z^{-j}\}:j\in \Z\}.
$$

\proof
If $\{z^{i_0}s\}$ is a basic set for some $i_0\in \mathbb{Z}\backslash \{0\}$,  applying Lemma \ref{1718c}, we obtain that $$
\{s,z^{2i_0}s\}=\{z^{i_0}s\}\{z^{-i_0},z^{i_0}\}
$$
is a basic set. This  yields that
$$
D=\{s,z^{2i_0}s\},
$$
and hence~$2i_0=i$. This yields that $i$ must be even and $\{z^{\frac{i}{2}}s\}$ is the unique basic set of size one consisting of nonidentity.
\medskip

Conversely, assume that $i$ is even. Then
$$
(s+z^is)(z^{\frac{i}{2}}+z^{-\frac{i}{2}})=z^{\frac{3i}{2}}s+z^{\frac{i}{2}}s+2z^{\frac{i}{2}}s
$$
belongs to $\cA$. Since $i\ne 0$, we obtain that $\{z^{\frac{i}{2}}s\}$ is a basic set by Corollary \ref{10131c}. Applying Lemma \ref{1718c}, one can see in this case that
$$
\cD=\{\{z^{j+\frac{i}{2}}s,z^{\frac{i}{2}-j}s\}:j\in \Z\}\}\cup \{\{z^j,z^{-j}\}:j\in \Z\}\}.
$$
 This completes the proof of Claim \textcolor[rgb]{0.00,0.00,1.00}{1}.
\eprf

\medskip

Now we may assume that~$i$ is odd. Then every basic set of nonidentity elements has size at least $2$ by Claim 1.  For any $j\in \Z\backslash \{0\}$, then both
$$
\{z^js,z^{i+j}s,z^{-j}s,z^{i-j}s\}=\{z^j,z^{-j}\}\{s,z^is\},
$$
and
$$
\{z^{-i-j}s,z^{i+j}s,z^{-j}s,z^{2i+j}s\}=\{s,z^is\}\{z^{i+j},z^{-i-j}\}
$$
are $\cA$-sets. It follows that
$$
\{z^{i+j}s,z^{-j}s\}=\{z^js,z^{i+j}s,z^{-j}s,z^{i-j}s\}\cap\{z^{-i-j}s,z^{i+j}s,z^{-j}s,z^{2i+j}\}
$$
is a basic set. As a consequence,
$$
\{z^{i+j}s,z^{-j}s\}\qaq \{z^js,z^{i-j}s\}
$$
are two basic sets.  Thus,
$$
\cD=\{\{z^js,z^{i-j}s\}:j
\in \Z\}\cup \{\{z^j,z^{-j}\}: j\in \Z\}.
$$

{\bf Subcase 2.3.}\, Suppose $|D|\ge 3$.  Then  there are $i\ne j\in\Z\backslash \{0\}$ such that $z^is$ and $z^js$ belong to $D$.
Observe that $C=\{z^i,z^{-i}\}$ is a basic set by the assumption. Since $s\in CD\cap D$, one can see that $D\subseteq CD$.

\medskip

{\bf Claim 2.}\, The structure constant $\lambda _{CDD}=1$.

\medskip

\proof
Otherwise, suppose the structure constant
$$
\lambda_{CDD}=|\{(x,y)\in C\times D:xy=z^is\}|\geq 2,
$$
then~$z^{2i}s\in D$.~ Observe that
$$
\lambda_{CDD}=|\{(x,y)\in C\times D:xy=z^{2i}s\}|\geq 2,
$$
which implies that~$z^{3i}s\in D$. It then follows that
that~$z^{ni}s\in D$ for all nonnegative integer $n$, which is a contradiction to the assumption that $D$ has finite support.  The claim then follows. \eprf

 \medskip

{\bf Claim 3.}\, The integers $i$ and $j$ satisfy  $j\neq2i,i\neq2j$.
\medskip

\proof Towards a contradiction, suppose that $j=2i$. Observe that
$$
z^i\cdot s=z^is,~z^{-i}\cdot z^{j}s=z^is,
$$
which yields that $\lambda _{CDD}\geq 2$,
which is a contradiction to Claim 2. So, $j\neq 2i$. Similarly,~we can conclude that~$i\neq 2j$ and the claim follows.

\eprf

Since $D\subseteq CD$  and $z^js\in D$,
one can see that
 \begin{equation}\label{2223a}
 z^{j-i}s\in D~~\mathrm{or}~~z^{j+i}s\in D.
\end{equation}
 Let $E=\{z^j,z^{-j}\}$.  Observe that~$s\in D\cap ED$, and hence $D\subseteq ED$.  Also $z^is\in D$. This implies that
 \begin{equation}\label{2223b}
 z^{i-j}s\in D~~\mathrm{or} ~~z^{j+i}s\in D.
\end{equation}
 If $z^{i+j}s\in D$, one can see that
 $$
 z^{i+j}s,s\in FD,
 $$
 where $F=\{z^{i+j},z^{-(i+j)}\}$.  So $D\subseteq FD$.
\medskip

 It follows that one of the following
 $$
 z^{-i}s,z^{-j}s, z^{i+2j}s, z^{2i+j}s
 $$
 belongs to $D$.
 \medskip

If $z^{-i}s\in D$, then
$$
z^i\cdot z^{-i}s=s,z^{-i}\cdot z^{i}s=s
$$
and thus
$$
\lambda_{CDD}=|\{(x,y)\in C\times D:xy=s\}|\geq 2,
$$
which is a contradiction to Claim 2.  So, $z^{-i}s$ does not belong to $D$. Similarly,  $z^{-j}s$ does not belong to $D$.

\medskip

Now, suppose that $z^{2i+j}s$ lies in $D$.  Then
$$
z^i\cdot z^js=z^{i+j}s,~z^{-i}\cdot z^{2i+j}s=z^{i+j}s.
$$
This yields that
$$
\lambda_{CDD}=|\{(x,y)\in C\times D:xy=z^{i+j}s\}|\geq 2.
$$
This is a contradiction to Claim 2. This means that $z^{2i+j}s$ does not lie in  $D$. Similarly $z^{i+2j}s$ does not lie in $D$. As a consequence,  $z^{i+j}s$ does not lie in $D$.

\medskip

Now, by statements \ref{2223a} and \ref{2223b} we conclude that $z^{i-j}s$ and $z^{j-i}s$ belong to $D$.  Let $K=\{z^{i-j},z^{j-i}\}$.~Then
$$
z^{i-j}\cdot z^{j-i}s=s,~z^{j-i}\cdot z^{i-j}s=s,
$$
and hence
$$\lambda_{KDD}=|\{(x,y)\in K\times D:xy=s\}|\geq 2.
$$
Observe that
$$
|\{(x,y)\in K\times D:xy=z^{i-j}s\}|=\lambda_{KDD}\geq 2,
$$
which implies that $z^{2(i-j)}\in D$. As in the proof of Claim 2, one can see that $z^{n(i-j)}s\in D$ for all nonnegative integers $n$. This is a contradiction to  $|D|<\infty$.~ Therefore neither $z^{i-j}s$ nor $z^{j-i}s$ lies in  $D$.  Consequently $|D|\leq 2$, i.e., subcase 2.3. does not occur.

\medskip

\eprf

Now
$$
\cD=\{\{z^js,z^{i-j}s\}:j
\in \Z\}\cup \{\{z^j,z^{-j}\}: j\in \Z\},
$$
where $D=\{s, z^is\}$ is the basic set containing $s$.  Let
$$
\varphi: ~G\longrightarrow G,~s\mapsto z^is,~z\mapsto z^{-1}.
$$
It is easy to see that $\varphi\in \mathrm{Aut}(G)$. Moreover,  $\cA={\F[\cZ\rtimes{\cZ}_2]}^{\langle\varphi\rangle}$.

\eprf
\medskip

{\bf Proof of Theorem \ref{2305b}}\, Now let
$
\mathcal{A}=\mathrm{Span}_{\mathbb{F}}\{~\underline{C}:C\in \mathcal{D}~\},
$
where
\begin{displaymath}
\cD=\{\{z^{i},z^{-i},z^is,z^{-i}s\}:2\nmid i,i\in \mathbb{Z}~\}\cup\{\{z^is^j, z^{-i}s^j\}:2\mid i,~i\in \mathbb{Z},j\in\{0,1\}~\}.
\end{displaymath}
Clearly, $\{e\}\in \mathcal{D}$ and $\mathcal{D}$ is closed to the Schur multiplication.
\medskip

 Let~$i_1,i_2\in \mathbb{Z}$~ with~$2\nmid i_1,2\nmid i_2$, then~$2|(i_1-i_2)$.~So,
\begin{displaymath}
\begin{split}
&(z^{i_1}+z^{-i_1}+z^{i_1}s+z^{-i_1}s)(z^{i_2}+z^{-i_2}+z^{i_2}s+z^{-i_2}s)\\
=&2(z^{i_1+i_2}+z^{-i_1-i_2})+2(z^{i_1+i_2}s+z^{-i_1-i_2}s)+\\&2(z^{i_1-i_2}+z^{i_2-i_1})+2(z^{i_1-i_2}s+z^{i_2-i_1}s)\in \mathcal{A}.
\end{split}
\end{displaymath}
Choose $i_3,i_4\in \mathbb{Z}$ with $2\nmid i_3,2\mid i_4$, then~$2\nmid(i_3-i_4)$. One can easily see that
\begin{displaymath}
\begin{split}
&(z^{i_3}+z^{-i_3}+z^{i_3}s+z^{-i_3}s)(z^{i_4}+z^{-i_4})\\
=&(z^{i_3+i_4}+z^{-i_3-i_4}+z^{i_3+i_4}s+z^{-i_3-i_4}s)+(z^{i_3-i_4}+z^{i_4-i_3}+z^{i_3-i_4}s+z^{i_4-i_3}s)\in \mathcal{A},
\end{split}
\end{displaymath}
\begin{displaymath}
\begin{split}
&(z^{i_4}+z^{-i_4})(z^{i_3}+z^{-i_3}+z^{i_3}s+z^{-i_3}s)\\
=&(z^{i_3+i_4}+z^{-i_3-i_4}+z^{i_3+i_4}s+z^{-i_3-i_4}s)+(z^{i_3-i_4}+z^{i_4-i_3}+z^{i_3-i_4}s+z^{i_4-i_3}s)\in \cA.
\end{split}
\end{displaymath}
Moreover,
\begin{displaymath}
\begin{split}
&(z^{i_3}+z^{-i_3}+z^{i_3}s+z^{-i_3}s)(z^{i_4}s+z^{-i_4}s)\\
=&(z^{i_3+i_4}+z^{-i_3-i_4}+z^{i_3+i_4}s+z^{-i_3-i_4}s)+(z^{i_3-i_4}+z^{i_4-i_3}+z^{i_3-i_4}s+z^{i_4-i_3}s)\in \cA,
\end{split}
\end{displaymath}
\begin{displaymath}
\begin{split}
&(z^{i_4}s+z^{-i_4}s)(z^{i_3}+z^{-i_3}+z^{i_3}s+z^{-i_3}s)\\
=&(z^{i_3+i_4}+z^{-i_3-i_4}+z^{i_3+i_4}s+z^{-i_3-i_4}s)+(z^{i_3-i_4}+z^{i_4-i_3}+z^{i_3-i_4}s+z^{i_4-i_3}s)\in \cA.
\end{split}
\end{displaymath}
Let $i_5,i_6\in \mathbb{Z}$ and $2\mid i_5,2\mid i_6$, $2\mid(i_5-i_6)$. Thus,
\begin{displaymath}
\begin{split}
&(z^{i_5}+z^{-i_5})(z^{i_6}+z^{-i_6})=(z^{i_5+i_6}+z^{-i_5-i_6})+(z^{i_5-i_6}+z^{i_6-i_5})\in \cA,\\
&(z^{i_5}+z^{-i_5})(z^{i_6}s+z^{-i_6}s)=(z^{i_5+i_6}s+z^{-i_5-i_6}s)+(z^{i_5-i_6}s+z^{i_6-i_5}s)\in \cA,\\
&(z^{i_5}s+z^{-i_5}s)(z^{i_6}+z^{-i_6})=(z^{i_5+i_6}s+z^{-i_5-i_6}s)+(z^{i_5-i_6}s+z^{i_6-i_5}s)\in \cA,\\
&(z^{i_5}s+z^{-i_5}s)(z^{i_6}s+z^{-i_6}s)=(z^{i_5+i_6}+z^{-i_5-i_6})+(z^{i_5-i_6}+z^{i_6-i_5})\in \cA.
\end{split}
\end{displaymath}
Hence, $\cA$ is a Schur ring over~$\cD_{\infty}$.
\medskip

Since there is no finite normal subgroup over $\cD_{\infty}$, this means that $\cA$ can not be the Schur rings of wedge product type. Obviously, there is no $\eta\in \mathrm{Aut}(\cD_{\infty})$ such that $\eta(z^i)=z^is$ for all $2\nmid i\in \mathbb{Z}$, which yields that $\cA$ is not an orbit Schur ring. So, $\cA$ is not traditional.

\medskip

\section{Acknowledgment}
The work of the first author is supported by Natural Science Foundation of China (No. 11971189, No. 12161035)


\begin{thebibliography}{99999}\small
	
\bibitem{BBHMT} N. Bastian, J. Brewer, M. Humphries, A. Misseldine, C. Thompson, On Schur rings over infinite
groups. Algebr. Represent. Theor. 23(3):493-511 (2020).

\bibitem{CP} G. Chen, I. Ponomarenko, \emph{Coherent Configurations}. Wuhan: Central China Normal Universit
Press (2019).

\bibitem{CH}
  G. Chen, J. W. He, Schur rings over $\cZ \times \cZ_3 $. Communications in algebra, 49(1): 4434-4446 (2021).


\bibitem{LM} K. H. Leung, S. L. Ma, The structure of Schur rings over cyclic groups. J. Pure Appl. Algebra 66(3):
287-302 (1990).

\bibitem{4}
 K. H. Leung, S. H. Man, On Schur rings over cyclic groups II. J. Algebra 183(2):273-285 (1996).




\bibitem{5} K. H. Leung, S. H. Man,  On Schur rings over cyclic groups. Isr. J. Math. 106(1):251-267 (1998).

\bibitem{Po}
I. Ponomarenko, \emph{Schur rings and algebraic combinatorics}. Lecture Notes. Wuhan (2015).



\bibitem{Sch}
I. Schur, \emph{Zur Theorie Der Einfach Transitiven Permutationsgruppen}, Vol. 118. Berlin: Sitzungsber.
Preuss. Akad. Wiss. Phy-Math Klasse, pp. 309-310 (1933).

\bibitem{Wie}
H. Wielandt, \emph{Zur theorie der einfach transitiven permutationsgruppen II}. Math. Z. 52(1):
384-393 (1950).

\bibitem{Wie2}
 H. Wielandt, \emph{Finite Permutation Groups}. New York, NY: Academic Press (1964).

\end{thebibliography}
\end{document}